# CORRELATION BETWEEN POLE LOCATION AND ASYMPTOTIC BEHAVIOR FOR PAINLEVÉ I SOLUTIONS

OVIDIU COSTIN

ABSTRACT. We extend the technique of asymptotic series matching to exponential asymptotics expansions (transseries) and show that the extension provides a method of finding singularities of solutions of nonlinear differential equations, using asymptotic information. This transasymptotic matching method is applied to Painlevé's first equation. The solutions of P1 that are bounded in some direction towards infinity can be expressed as series of functions, obtained by generalized Borel summation of formal transseries solutions; the series converge in a neighborhood of infinity. We prove (under certain restrictions) that the boundary of the region of convergence contains actual poles of the associated solution. As a consequence, the position of these exterior poles is derived from asymptotic data. In particular, we prove that the location of the outermost pole $x_p(C)$ on $\mathbb{R}^+$ of a solution is monotonic in a parameter $C$ describing its asymptotics on antistokes lines[1], and obtain rigorous bounds for $x_p(C)$. We also derive the behavior of $x_p(C)$ for large $C \in \mathbb{C}$. The appendix gives a detailed classical proof that the only singularities of solutions of P1 are poles.

## 1. INTRODUCTION AND NOTATIONS

We consider the solutions of Painlevé's first equation, P1

$$(1.1) \qquad y'' = 6y^2 + z$$

and study the connection between their asymptotic behavior and the location of their poles. Classical asymptotic analysis and isomonodromy theory [1] have not yielded so far this type of information.

The analysis and proofs in the present paper are mostly self-contained, but the method that we use, as well as the nature of results obtained, are much easier to interpret from the more general perspective of exponential asymptotics. We start by briefly presenting some concepts and results in this area which are relevant to the problem under discussion.

Jean Écalle has introduced and developed the theory of analyzable functions and tools of generalized resummation of *transseries* (generalized exponential-power series) as a natural and very comprehensive generalization of the theory of analytic functions [2, 3, 4]. The grand picture is the construction of a class of functions, closed

---

Research at MSRI was supported, in part, by NSF grant DMS-9704968.
[1] In some texts called Stokes lines: the lines where the exponentials in the asymptotic expansion are purely oscillatory.

under all imaginable algebraic and analytic operations, which are into an isomorphic correspondence with the space of transseries. In a sense, the space of transseries represents the formal closure of power series under all operations. A wide class of analytic problems admit formal solutions in the space of transseries (obtained in a more or less straightforward way). Transseries solutions contain ample analytic information and in many cases (for instance generic linear or nonlinear differential systems [8]) the actual solutions of the problem can be reconstructed from their transseries by modern resummation methods.

1.1. **Notes on asymptotics beyond all orders for P1.** Proofs of the various technical statements in this section can be found in [8], in the more general context of analytic rank-one differential systems.

The change of variables

$$(1.2) \qquad i\sqrt{6z^{-1}}y(z) = -1 + \frac{4}{25x^2} - h(x); \quad 30x = (-24z)^{5/4}$$

brings (1.1) to the (normal) form

$$(1.3) \qquad h'' + \frac{1}{x}h' - h - \frac{1}{2}h^2 - \frac{392}{625x^4} = 0$$

There is a one-parameter family of solutions of (1.3) that tend to zero as $x \to +\infty$, $x \in \mathbb{R}$, and they have a common asymptotic series (cf. also § 3):

$$(1.4) \qquad h(x) \sim \tilde{h}_0(x) = \sum_{k=4}^{\infty} c_k x^{-k}$$

The classical expansion (1.4) does not distinguish different solutions in this family. However there exists a one complex parameter family of transseries solutions, of which $\tilde{h}_0$ is only a special instance:

$$(1.5) \qquad \tilde{h}(x, C) = \tilde{h}_0(x) + \sum_{k=1}^{\infty} C^k x^{-k/2} e^{-kx} \tilde{h}_k(x) \quad (C \in \mathbb{C})$$

In (1.5), $\tilde{h}_k$ are formal power series, $\tilde{h}_k = \sum_{m=0}^{\infty} c_{km} x^{-m}$, that can be computed algorithmically by substituting (1.5) into (1.3) and identifying the coefficients. With the choice $c_{10} = 1$, which is possible for (1.3), the series $\tilde{h}_k$ do not depend on $C$. Every $\tilde{h}_k$ is factorially divergent. There is nevertheless a very close connection between (1.5) and the family of actual solutions of (1.3) which decrease on $\mathbb{R}$ for $x \to +\infty$: any such solution can be reconstructed by generalized Borel resummation of $\tilde{h}(x, C)$ for some $C$. Denote by $\mathcal{LB}$ the generalized Borel resummation operator (as defined in [8]). There exists $x_0 > 0$ (independent of $k$) so that the functions $h_k = \mathcal{LB}\tilde{h}_k$ are analytic in the half plane $H_{x_0} = \{x : \Re(x) > x_0\}$. In addition, $h_k \sim \tilde{h}_k$ as $x \to e^{i\theta}\infty$ in $H_{x_0}$, $\theta \in (-\pi/2, \pi/2)$. For any $C \in \mathbb{C}$, the function series

$$(1.6) \qquad h(x;C) = \mathcal{LB}\tilde{h}_0(x) + \sum_{k=1}^{\infty} C^k x^{-k/2} e^{-kx} \mathcal{LB}\tilde{h}_k(x)$$

converges for $\Re(x) > x_1 = x_1(C)$ to a solution of (1.3); conversely, any solution $h$ which vanishes for real $x \to +\infty$ is necessarily of the form (1.6) for a unique $C$. $\mathcal{LB}$ commutes with the algebraic operations, complex conjugation (thus a series with real coefficients sums to a real function), with differentiation, and has good continuity properties. Along directions in $\mathbb{C}$ other than $\mathbb{R}^+$ similar results hold.

In addition, the constant $C$ in (1.6) is related by a simple formula to the large $|x|$ behavior of $h(x;C)$ on the antistokes lines $\pm i\mathbb{R}^+$ (cf. (2.8)). We show that for any $C > 0$, the boundary of the domain where the series (1.6) converges contains an actual pole of $h(x;C)$, much as the boundary of the convergence disk of the Taylor series of an analytic function must contain a singularity.

1.2. **Heuristic discussion of transasymptotic matching.** The value of $C$ is related to the position of the outermost singularities of solutions of P1. The following heuristic argument illustrates the reason.

For simplicity, we assume that $C$ is large. When $x$ is very large, the exponentials in (1.5) are smaller than the leading term in $\tilde{h}_0$, and the behavior of $h(x;C)$ is roughly $c_{04}x^{-4}$.

For smaller values of $x$, such that $1 \gg |Ce^{-x}x^{-1/2}| \gg |x|^{-1}$, the relative magnitude of the terms of (1.5) changes: the largest term becomes $Ce^{-x}x^{-1/2}$ (instead of $c_{04}x^{-4}$), the next order of magnitude is the term $C^2 e^{-2x} x^{-1} c_{20}$ (not $Ce^{-x}x^{-3/2}c_{11}$) and so on. In this intermediate regime, the terms of (1.5) need to be reshuffled accordingly, and $h(x;C)$ is represented by

$$(1.7) \qquad \tilde{h}(x;C) = \sum_{m=0}^{\infty} \frac{1}{x^m} \left( \sum_{k=0}^{\infty} \frac{c_{km}}{(e^x x^{1/2}/C)^k} \right) := \sum_{m=0}^{\infty} \frac{G_m(e^x x^{1/2}/C)}{x^m}$$

It turns out that, for $x$ not too small, the inner sums converge (while, obviously, the expansion $\tilde{h}$ as a whole is still divergent), and $G_m(t)$ are analytic. In fact, for (1.3), $G_m(t)$ are rational functions and $G_0$ has a particularly simple expression, $G_0(t) = t(t - 1/12)^{-2}$ (see the proof of Lemma 8). Moreover, for any $m \geq 0$ the function $G_m(t)$ has only one pole in the right half plane, at $t = 1/12$.

The representation (1.7) is normally accurate as long as no Stokes phenomena occur and the following asymptoticity condition holds

$$(1.8) \qquad \left| G_{m+1}(e^x x^{1/2}/C) x^{-m-1} \right| \ll \left| G_m(e^x x^{1/2}/C) x^{-m} \right|$$

A formal calculation shows that (1.8) holds for $x$ in a region $\mathcal{R}$ that starts well to the *left* of the array of poles of $g(x) = G_0(e^x x^{1/2}/C)$, located at $x_N$, where

$$(1.9) \qquad x_N + \frac{1}{2}\ln x_N = \ln(C/12) + 2N\pi i; \ N \in \mathbb{Z}$$

and extends all the way (small neighborhoods of the poles themselves excluded) to $\Re(x) \to \infty$. The region $\mathcal{R}$ is thus substantially wider than the region where (1.4) is valid.

Accepting for the moment that (1.7) is the asymptotic expansion of some solution $h$ in $\mathcal{R}$, (1.7) can be used to estimate $\oint(\int h)$ on a circuit around a pole $x_N$. The calculation gives $\oint(\int h) = -12\pi i + O(1/\ln C)$ for large $|C|$, showing that $h$ is singular at some point $x_p(C; N) = x_N + o(1)$. $x_p(C; N)$ can be estimated to all orders in $1/\ln C$ from (1.7).

On the other hand, if $x$ is large, and away from $x_N$, $h(x; C)$ must be bounded (by (1.7)). The location of the poles of $h$ in $\mathcal{R}$ should thus be $x_p(C; N) = x_N + o(1)$.

This method gives information not only about the location, but also about the nature of the singularities of $h(x; C)$ in $\mathcal{R}$. A straightforward calculation based on the Cauchy integral argument explained above indicates, to any order in $1/\ln C$, that $x_p(C; N)$ are double poles.

The fact that P1 is an integrable equation plays no role, except that $G_m$ can be computed in closed form.

We note that in the matching process *all the terms* of the transseries (1.5) are necessary, so that the matching is *trans*asymptotic. The procedure makes the connection between the behavior of $h$ in the pole-free region near infinity, and its singularities for large $x$ in $\mathbb{C}$.

In the next sections these conclusions are made rigorous, and further results are proved on the connection between $C$ and the position of $x_p(C; 0)$, for $C$ not necessarily large. The proofs yield global estimates for (1.6) and (1.7).

We can also interpret (1.7) from the perspective of multiscale analysis. Apparently an expansion of the form $\sum c_m(x) x^{-m}$ (where the coefficient of $x^{-m}$ may depend on $x$) is subject to substantial arbitrariness. Not so, however, if $c_m(x)$ are required to be sufficiently regular and the expansion to be valid over a wide enough range of $x$. This is a crucial point in multiscale analysis. We note that (1.5) determines the asymptotic scales $x^{-1}$ and $Ce^{-x}x^{1/2}$. Once the scales are found, the terms in (1.7) are more conveniently determined by inserting (1.7) in (1.3) and *identifying the coefficients of $x^{-m}$*.

The proper choice of the asymptotic scales is essential. For instance looking for a representation of the more general form $\tilde{h}_a = \sum x^{-m} G_m(Ce^{-x}x^a)$, only for $a = 1/2$ will the expansion be asymptotic (even formally) for large $x$. For any other value of $a$, $G_m(Ce^{-x}x^a) = O(x^m)$ instead of $O(1)$, and then $\tilde{h}_a$ is rather meaningless.

## 2. Main results

*Notations.* By Lemma 2 in § 3 there is a one (real) parameter family of solutions of (1.3) which are real-valued for $x > const.$ on $\mathbb{R}^+$ and such that $h_0(x) \to 0$ as $x \to +\infty$. Their behavior for large $x$ is $h_0 \sim c_4 x^{-4} + c_6 x^{-6}(1 + o(1))$ where $c_4 = -392/625$ and $c_6 < 0$. We choose $h_0$ to be one such solution. (The choice $h_0 = \mathcal{LB}\tilde{h}_0$ has however important additional properties). We let $A \geq 9/4$ be such that, for $x > A$,

$$(2.1) \qquad h_0(x) \in \left(-x^{-4}, c_4 x^{-4}\right)$$

**Proposition 1.** *(i) Any solution $h$ of (1.3) such that $|h(x)| < 1$ for large positive $x$ is given by*

$$h(x;C) = \sum_{k=0}^{\infty} C^k x^{-k/2} e^{-kx} h_k(x) \tag{2.2}$$

*for some $C \in \mathbb{C}$, where the functions $h_k(x)$ are determined recursively as solutions of linear differential equations (cf. (3.2)) and do not depend on $C$.*

*The series (2.2) converges for all $x > x_0(C) = \max\{A, x_1\}$ where $x_1 + \frac{1}{2}\ln x_1 = \ln(|C|/12)$.*

*(ii) The functions $h_k$ satisfy the inequalities*

$$0 < h_k(x) < \lim_{x \to \infty} h_k(x) = \frac{k}{12^{k-1}} \quad (k \in \mathbb{N} \text{ and } x > A) \tag{2.3}$$

*Furthermore,*

$$\overline{\lim_{k \to \infty}} h_k(x)^{-1/k} = 12 + O(x^{-1}) \text{ (as } x \to +\infty)$$

*(iii) (Outermost pole of $h(x;C)$ on $\mathbb{R}^+$.) If $C > C_0 = 12 \min_{t \geq A} \sqrt{t} e^t (1 - \frac{9}{4t})^{-1}$, then $h(x;C)$ has a double pole at the point $x_p = x_p(C; 0)$, the unique solution of*

$$x^{1/2} e^x \overline{\lim_{k \to \infty}} h_k(x)^{-1/k} = C \tag{2.4}$$

*for $x > A$. The function $h(x;C)$ has no other singularities on $[x_p, \infty)$.*

*(iv) (Monotonicity and bounds for $x_p$.) For $C > C_0$, $x_p(C; 0)$ is increasing in $C$ and satisfies the inequalities*

$$x_p + \frac{1}{2}\ln x_p \leq \ln(C/12) \leq x_p + \frac{1}{2}\ln x_p - \ln\left(1 - \frac{9}{4x_p}\right) \tag{2.5}$$

*Thus, $x_p$ has the asymptotic behavior,*

$$x_p = \ln(C/12) - \frac{1}{2}\ln(\ln(C/12)) + O(1/\ln C) \ (C \to +\infty) \tag{2.6}$$

*(v) (Poles of $h(x;C)$ in $\mathbb{C}$, for large $C$ complex.) Let $N \in \mathbb{N}$ and $\epsilon > 0$. As $|C| \to \infty$, $h(x;C)$ has double poles at the points:*

$$x_p(C; k) = \ln(C/12) - \frac{1}{2}\ln(\ln(C/12)) + 2k\pi i + o(1) \quad (k \in \mathbb{Z} \cap [-N, N]) \tag{2.7}$$

*and no poles in the region $\{x : |\Im(x)| \leq N, \ \Re(x) > \Re(x_p(C; 0)) + \epsilon\}$.*

*(vi) (Connection to behavior on antistokes lines.) For the special choice $h_0 = \mathcal{LB}\tilde{h}_0$ we have $h_k = \mathcal{LB}\tilde{h}_k$, and*

$$h(x;C) \sim (C \pm \frac{1}{2}S) x^{-1/2} e^{-x} \text{ as } x \to \pm i\infty \tag{2.8}$$

where $C \in \mathbb{C}$ is the same as in (2.2) and $S$ is the Stokes constant of the equation (1.3).

**Notes**. (1) The property that the solutions which vanish in some direction can be expressed, as in (2.2), as a power series in the constant $C$ of (2.8) is very general, cf. [8], and explains why $C$ is also a natural parameter for the location of exterior poles. Integrability of P1 plays no role here. A formula of type (2.4) is expected to be quite general and is not a connection formula stricto sensu, because it involves limits. For P1 however, there are indications that these limits can be calculated in closed form.

(2) For the choice $h_0 = \mathcal{LB}\tilde{h}_0$ it can be shown, using the methods in [9], that (2.1) holds for $x > 4/3$.

## 3. Proofs and further results

Before proving the main proposition we describe some asymptotic properties of the solutions of (1.3) which are bounded as $x \to +\infty$.

**Lemma 2.** *(i) If $h$ is a solution of (1.3) and $|h(x)| < 1$ for large positive $x$, then*

$$(3.1) \qquad h \sim \tilde{h}_0 = \sum_{k=4}^{\infty} c_k x^{-k} \qquad (x \to e^{i\theta}\infty \text{ with } \theta \in \left(-\frac{\pi}{2}, \frac{\pi}{2}\right))$$

*(the first few terms of $\tilde{h}_0$ are $-\frac{392}{625}x^{-4} - \frac{6272}{625}x^{-6} - \frac{141196832}{390625}x^{-8} - ...$).*
*(ii) There is a one (complex/real) parameter of such solutions (complex/real valued on $\mathbb{R}^+$); if $h_1$ and $h_2$ satisfy (3.1) then $h_1 - h_2 = Ce^{-x}x^{-1/2} + o(e^{-x}x^{-1/2})$, for some $C \in \mathbb{C}$, as $x \to \infty$. If in addition $C = 0$ then $h_1 = h_2$.*

*Proof of Lemma 2.* (i) is standard (see, e.g. [6]); it is also easy to reprove in our case, since (1.3) can be rewritten as

$$h'' + \frac{1}{x}h' - \left(1 + \frac{1}{4x^2}\right)h = \frac{1}{2}h^2 - \frac{1}{4x^2}h + \frac{392}{625x^4} := g(x)$$

which, in integral form

$$h(x) = x^{-1/2}\left(e^x \int_a^x e^{-t}t^{1/2}g(t)\mathrm{d}t + e^{-x}\int_\infty^x e^{-t}t^{1/2}g(t)\mathrm{d}t\right)$$

is a contraction in the unit ball in the sup norm, for large $x$. If $a \in \mathbb{R}$ then $h$ is real-valued.

Part (ii) is proven without difficulty in essentially the same way, by writing an integral contractive equation for $\delta = h_2 - h_1$ (see [8] for details under more general assumptions). □

3.1. **Properties of (2.2) on $\mathbb{R}^+$.** Let $h_0$ and $A$ be as in § 2. Formal substitution of (2.2) in (1.3) yields

$$(3.2) \qquad H_k'' - \left(1 + h_0(x) - \frac{1}{4x^2}\right) H_k = \frac{\sqrt{x}}{2} e^{-kx} \sum_{0 < j < k} h_j h_{k-j} := R_k(x) \ (k \geq 1)$$

where $H_k(x) = x^{-(k-1)/2} e^{-kx} h_k(x)$.

*Step 1.* Study of the homogeneous part of (3.2) (same as (3.2) for $k = 1$):

$$(3.3) \qquad y'' - \left(1 + h_0(x) - \frac{1}{4x^2}\right) y = 0$$

**Lemma 3.** *(i) Eq. (3.3) has two solutions $y_+, y_-$ such that: $y_\pm(x) > 0$ for $x > A$, $y_\pm(x) = e^{\pm x}(1 + O(x^{-1}))$ as $x \to +\infty$ and their Wronskian $(y_- y_+' - y_-' y_+)$ is $W = 2$.*
*If $y(x)$ solves (3.3) and $y(x) = O(e^{-x})$ as $x \to \infty$ then $y = const. y_-$.*
*(ii) The function $u(t, s) = y_+(t) y_-(s) - y_-(t) y_+(s)$ is positive for $t > s > A$.*

*Proof.* (i) Except for the explicit region of positivity of $y_\pm$, this is standard asymptotic theory (see e.g. [6]; note also that $W = const$ and $W \to 2$ for $x \to \infty$). We prove the statement only for $y_-$; the argument for $y_+$ is similar. Associating to (3.3) the integral equation

$$(3.4) \qquad f(x) = 1 - \int_x^\infty e^{x-s} \sinh(s - x) \left(\frac{1}{4s^2} - h_0\right) f(s) ds$$

where $f(x) = e^x y(x)$, we define

$$(3.5) \qquad J(\Psi)(x) = \int_x^\infty e^{x-s} \sinh(s - x) \left(\frac{1}{4s^2} - h_0\right) \Psi(s) ds$$

on $L^\infty(x_0, \infty)$, for $x_0 > A$.

**Remark 4.** $\|J\|_{L^\infty \mapsto L^\infty} < \frac{1}{8x_0}\left(1 + \frac{1}{3x_0^2}\right)$ *and $J$ is positivity preserving. Thus (3.4) has a unique solution $f_1$ in $L^\infty(A, \infty)$; moreover, $0 < f_1 < 1$ and $f_1(x) \sim 1 + \sum_{j=1}^\infty (-1)^j |c_j| x^{-j}$ $(x \to +\infty)$.*

*Proof.* The first part is straightforward. Now, $\|f_1\|_\infty \leq (1 - \|J\|)^{-1} \|1\|_\infty$ thus $|f_1| \leq \frac{6}{5}$ implying $|f_1 - 1| \leq \frac{6}{5} J(1) \leq \frac{1}{5}$, in particular $f_1 > 0$. Since $J$ preserves positivity, we have $f_1 = 1 - J f_1 < 1$. The asymptotic expansion of $f_1$ follows from the representation $f_1 = \left(\sum_{j=0}^\infty (-J)^j\right) 1$. □

For the proof of Lemma 3 (ii) note that $u(t, s)$ is a solution of $u_{tt} = (1 + h_0(t) - \frac{1}{4} t^{-2}) u$ with initial condition $u(s, s) = 0$; $u_t(s, s) = 2$. Since $1 + h_0(t) - \frac{1}{4} t^{-2} > 0$ for $t > A$, we have $u(t, s) > 0$ for all $t > s$. □

**Lemma 5.** *There is a unique solution of (3.2) with $H_1 = y_-$ such that for $k > 0$, $H_k(x) = O(e^{-kx})$ as $x \to \infty$. Namely, the $H_k(x)$ are inductively given by*

$$(3.6) \qquad H_k(x) = \frac{1}{2}\int_x^\infty u(t,x) R_k(t)\mathrm{d}t, \ k \geq 2$$

*In addition $H_k(x) > 0$ for $x > A$.*

*Proof.* Immediate induction using Lemma 3. Uniqueness follows from the behavior for large $x$ of the solutions of the homogeneous equation, $y_\pm$ (Lemma 3). □

*Step 2. Upper bounds for $H_k$.* We have, by (3.2) and Lemma 5

$$(3.7) \quad H_k(x) = \int_x^\infty \sinh(t-x) R_k(t) dt - \int_x^\infty \sinh(t-x)\left(\frac{1}{4t^2} - h_0(t)\right) H_k(t) dt$$

$$\leq \int_x^\infty \sinh(t-x) e^{-kt} R_k(t) dt$$

From Lemma 3 and Remark 4 we get $H_1(x)e^x \in (0,1)$ on $(A,\infty)$. Assuming by induction that $H_m(x) \leq \frac{m}{12^{m-1}} x^{-(k-1)/2} e^{-mx}$ if $1 \leq m < k$, we get

$$(3.8) \quad H_k(x) \leq \frac{1}{2}\int_x^\infty \sinh(t-x) t^{-k/2+1/2} e^{-kt} \sum_{0<j<k} \frac{j(k-j)}{12^{k-1}} dt \leq$$

$$\frac{(k^3-k)}{12^k} x^{-(k-1)/2} \int_x^\infty \sinh(t-x) e^{-kt} dt = \frac{k}{12^{k-1}} x^{-(k-1)/2} e^{-kx}$$

□

*Step 3. Monotonicity*

**Lemma 6.** *The functions $H_k(x)$ and $F_k(x) = x^{-1/2} e^{-kx} h_k(x)$, $k \geq 1$, are decreasing for $x > A$.*

*Proof.* The functions $F_k$ satisfy the equations

$$(3.9) \qquad F_k'' = (1 + h_0) F_k - \frac{1}{x} F_k' + \frac{1}{2} \sum_{0<j<k} F_j F_{k-j} \qquad (k \geq 1)$$

By (3.8) and Lemma 5, the $F_k$ are positive and vanish as $x \to \infty$. On the other hand, $F_k$ can have no extremum for $x > A$ (otherwise, by (3.9) all extrema would be *minima*), thus are decreasing in $x$. For $H_k$ the proof is very similar. □

*Step 4. Lower bounds for $H_k$.*

**Lemma 7.** *For $x > A$ we have*

$$(3.10) \qquad H_m(x) \geq \left(1 - \frac{1}{8x}\right)\left(1 - \frac{9}{4x}\right)^{(m-1)/2} \frac{m}{12^{m-1}} x^{-(k-1)/2} e^{-mx}$$

From Remark 4 we get $f = 1 - Jf \geq 1 - J1 \geq 1 - (1 - (8x)^{-1})^{-1}$ for $x > 5/3$. Thus $H_1(x) \geq (1 - (8x)^{-1})$. Since the functions $H_k$ are decreasing, we have

$$(3.11)\quad H_k(x) \geq \int_x^\infty dt \sinh(t-x) R_k(t) - H_k(x) \int_x^\infty dt \sinh(t-x) \left(\frac{1}{4t^2} - h_0(t)\right)$$

Thus

$$(3.12)\quad H_k(x) \geq \left(1 - \frac{1}{8x}\right) \int_x^\infty \sinh(t-x) R_k(t) dt$$

From this point on, the proof of Lemma 7 follows the same track as Step 2, noting that $(1 - t^{-1})^{(k-1)/2}$ is increasing in $t$ and that $\int_x^\infty \sinh(t-x) t^{-(k-1)/2} e^{-kt} \geq (1 - x^{-1})(k^2 - 1)^{-1} x^{-(k-1)/2} e^{-kx}$. □

*Step 5. Completing the proof of Proposition 1 (i) and (ii).* With $C \in \mathbb{C}$ given, (2.2) converges for $x > x_0(C)$, by (3.8). Furthermore, since the $H_k$ are decreasing, if (2.2) converges at $x = x_1$, then it converges uniformly on $[x_1, \infty)$ (together with its derivatives, in view of the equality in (3.7)) and thus is a solution of (1.3). By (3.8) and (3.10), we have $h(x; C) - h_0(x) = Ce^{-x} x^{-1/2}(1 + o(1))$ for large $x$. Therefore, by Lemma 2, if $g$ is solution of (1.3) which vanishes as $x \to \infty$, then $g = h(x; C)$ for some $C \in \mathbb{C}$. □

*Step 6. Divergence of the series (2.2) implies singularities of the function $h(x; C)$.* Let $C > C_0$ and $x_{div} := \inf\{x | \sum_{k=1}^\infty C^k F_k(x) < \infty\}$. We have just shown that (2.2) solves (1.3) for $x > x_{div}$. By the definition of $C_0$ and (3.7) we have $x_{div} > A$. We now want to show that the solution $h(x; C)$ is nonanalytic at $x_{div}$. Let $S_N(x) = \sum_{k=1}^N C^k F_k(x)$.

**Lemma 8.** *Assume $h(x; C)$ is a $C^2(x_1, \infty)$ solution of (1.3), and $x_1 \in (A, x_{div}]$.*
*(i) The function $\delta_N(x) = h(x; C) - h_0 - S_N(x)$ is positive for all $x > x_1$ and $N \in \mathbb{N}^+$.*
*(ii) $x_1 = x_{div} = x_p(C; 0)$.*

*Proof.* (i) We have, by Proposition 1 (ii) $\delta_N(x) = \sum_{j>N} C^j F_j(x) > 0$ for $x > x_{div}$. Now, since the function $\delta_N$ satisfies the equation

$$(3.13)\quad \delta_N'' = \frac{1}{2}\delta_N^2 + (1 + S_N + h_0)\delta_N - \frac{\delta_N'}{x} + \frac{1}{2} \sum_{\substack{i,j \in (1, N-1) \\ i+j > N}} F_i F_j$$

we see, as in Lemma 6, that $\delta_N$ is decreasing, therefore positive, for all $x > x_1$.
(ii) By (i), $h(x; C)$ is unbounded at $x_{div}$. But the only singularities of $h(x; C)$ for $x \neq 0$ are double poles (see the appendix). □

*Step 7. Concluding the proof of Proposition 1, (iii) and (iv).* The monotonicity of $x_p(C; 0)$ in $C$ is now immediate, while (2.5) follows from Lemma 8, (3.8) and (3.10). □

**Note.** Formula (2.4) in itself is not a practical numerical way to calculate $x_p(C)$, given a particular $C$. An efficient (and rigorously controllable) method is to use the

procedure in this section, to obtain successively better upper and lower bounds for $H_k$.

## 3.2. Correspondence between (1.6), (2.2) and (2.8).

**Lemma 9.** *If $h_0 = \mathcal{LB}\tilde{h}_0$ then $h_k = \mathcal{LB}\tilde{h}_k$ for all $k \in \mathbb{N}$. Thus the representations (1.6) and (2.2) coincide. Formula (2.8) holds.*

*Proof of Lemma 9 and of Proposition 1, (vi).* We let $g_{1,2}(x,C) = h(x;C)$ as given by (1.6) and (2.2), respectively. In [8], formula (2.8) for $g_1(x,C)$ and the estimate $g_1(x,C) - h_0(x) = Ce^{-x}x^{-1/2}(1+o(1))$ are proved in a general setting. Lemma 2 then implies, as in Step 5 § 3.1, that $g_1(x;C) = g_2(x;C)$. □

**Note.** Using Lemma 9 it follows from [8] that (2.2) converges in fact for all large $x$ with $\arg(x) \in (-\pi/2, \pi/2)$.

## 3.3. Uniform representation for large $x$.

We are now interested in the region of validity of (1.7) and deriving the behavior of $x_p(C)$ for large $C \in \mathbb{C}$. While the results in the previous sections used specific properties of (1.3) (positivity of $h_k$), the technique of transasymptotic matching introduced in §1.2 is expected to work for more general nonlinear systems. We will only prove (1.7) to leading order, sufficient to obtain (2.7); cf. also Remark 12. The expansion (1.7) is studied in a region of $x$ containing points near $x_N$ (see (1.9)) and extending to infinity. Let $\epsilon > 0$ be small, $M \in \mathbb{R}^+$ and $s = s(x;C) = 12e^x x^{1/2}/C$. Define $\mathcal{R} = \{x : |\Im(x)| \leq M\,|s-1| \geq \epsilon$ and $\Re(s) > 1 - \epsilon\}$.

**Lemma 10.** *Let $C_1 \in \mathbb{R}^+$ be large. For any $C$ with $|C| > C_1$ we have*

$$(3.14) \qquad h(x;C) = \frac{12s}{(s-1)^2} + \psi(x;C)$$

*where $|\psi(x;C)| \leq K|\ln C|^{-1}$ in $\mathrm{Clos}(\mathcal{R})$, for some $K$ not depending on $C$.*

*Step 1: Proof of (3.14) for $x \in \mathbb{R}^+$.* From (2.2) and Lemma 7 we obtain

$$(3.15) \quad h(x;C) = h_0 + \sum_{k=1}^{\infty} 12k s^{-k} + \sum_{k=1}^{\infty} 12k s^{-k}\left[(1 + r_k(x)x^{-1})^k - 1\right] =$$

$$\frac{12s}{(s-1)^2} + \sum_{k=1}^{\lfloor x \rfloor} \frac{12k^2 c_k(x)}{x} s^{-k} + \sum_{\lfloor x \rfloor + 1}^{\infty} 12k\left[(1 + r_k(x)x^{-1})^k - 1\right] + O(x^{-4}) =$$

$$\frac{12s}{(s-1)^2} + \frac{E_1(x)}{\ln(C)} \quad \text{with } |E_1(x)| \leq const. \ (x \in \mathbb{R}^+ \cap \mathcal{R})$$

The estimates in (3.15) are obtained in the following way: by Lemma 7 we have $|r_k(x)| \leq 9/4$; thus, $|c_k(x)| = |(1 + r_k/x)^k - 1| \leq (e^{9/4} - 1)kx^{-1}$ for $k \leq x$; the tail $\sum_{\lfloor x \rfloor + 1}^{\infty}(\cdot)$ is exponentially small, while for $x \in \mathcal{R}$ we have $h_0 = O(x^{-4})$ and $|x/\ln C| \leq 2$ as $C \to \infty$. Similarly,

$$(3.16) \qquad \left|\frac{d}{dx}\left(h(x,C) - 12s(s-1)^{-2}\right)\right| \le \frac{const}{\ln(C)} \text{ for } x \in \mathbb{R}^+ \cap \mathcal{R})$$

*Step 2. Extension to $\mathcal{R}$.* Let $\psi(x) = h(x;C) - 12s/(s-1)^2$. It is convenient to change variables to $v = x + \frac{1}{2}\ln x - \ln(C/12)$ and define $\mathcal{R}_v = \{v|x \in \mathcal{R}\}$. Then $\boldsymbol{\psi} = (\psi, \psi')$ satisfies the system

$$(3.17) \qquad \frac{d}{dv}\begin{pmatrix}\psi_1\\\psi_2\end{pmatrix} = \begin{pmatrix}0 & 1\\F(v) & 0\end{pmatrix}\begin{pmatrix}\psi_1\\\psi_2\end{pmatrix} + \begin{pmatrix}0\\\left(1+\frac{1}{2x}\right)^{-2}r\end{pmatrix} := A\boldsymbol{\psi} + \mathbf{R}$$

where $F(v) = (10e^v + e^{2v} + 1)(e^v - 1)^{-2}$ and

$$r = \frac{\psi_1^2}{2x^{1/2}} + \frac{(x+\frac{1}{4})\psi_1}{(x+\frac{1}{2})^2} + \frac{2\psi_2 - \psi_1}{4x^2} - 3\frac{e^{3v} + 4e^{2v} + e^v}{(e^v-1)^4 x^{3/2}} + \frac{24e^v(2e^v+1)}{x^{1/2}(e^v-1)^4} - \frac{c_4}{x^{7/2}}$$

Let $U(t;v)$ be the fundamental matrix of $U' = AU$ (the solution with $U(t;t) = I$, the identity matrix). [2]

**Remark 11.** *$U(t;v)$ and $U^{-1}(t;v)$ are uniformly bounded on*

$$\mathcal{R}_{M;t} = \{(t,v) \in \mathcal{R}_v^2 | \ |t-v| < M\}$$

*Proof.* $U$ is uniformly bounded since $F$ is uniformly continuous on $\mathcal{R}_v$. Since $\text{Tr}(A) = 0$ we have $\det(U) = const = 1$, so $U^{-1}$ is bounded as well. □

Let $v_0 \in \mathcal{R}_v \cap \mathbb{R}$. By (3.15) and (3.16) we have $|\boldsymbol{\psi}(v_0)| \le const./\ln C$. We write (3.17) in integral form

$$(3.18) \qquad \boldsymbol{\psi}(v) = U(v_0;v)\boldsymbol{\psi}(v_0) + U(v_0;v)\int_{v_0}^v U^{-1}(v_0;t)\mathbf{R}(t)dt$$

The integral equation (3.18) is contractive in $L^\infty(\mathcal{R}_{M;v_0})$ in a ball of radius $K/\ln C$ for some $K$ which is, by (3.15) (3.16) and Remark 11, independent of $v_0$ (note that all terms on the r.h.s. of (3.18) are in fact $O(1/\ln C)$ if $\|\boldsymbol{\psi}\| = 1$). Lemma 10 is proven. □

*Step 3: Pole location for large $C \in \mathbb{C}$. Proof of Proposition 1, (v).* Let

$$\mathcal{K} = \{x : |s-1| = \epsilon \text{ and } |v - 2N\pi i| < 2\pi\}$$

If $\epsilon$ is small enough, then for large $C$ the curve $\mathcal{K}$ is a simple closed curve in $\text{Clos}(\mathcal{R})$. We have by Lemma 10 $\oint_\mathcal{K}\left(\int h(\cdot;C)\right) = -12\pi i + o(1)$ thus $h$ is not analytic in the interior of $\mathcal{K}$. The only singularities of $h(x;C)$ for $x \ne 0$ are double poles, cf. the appendix. □

---

[2] Incidentally, the equation $U' = AU$ can be solved explicitly: for the associated second order equation the general solution is $\psi = C_1 t(t+1)(t-1)^{-3} + C_2(t^{-1}P_5 + 60tv(t+1))(t-1)^{-3}$, where $t = e^v$ and $P_5$ is a quintic polynomial in $t$.

**Remark 12.** *Using the asymptotic behavior of the solution of $U' = AU$ and more careful estimates, it is possible to prove contractivity and thus (3.14) in a sector of opening less than $\pi$ centered on $\mathbb{R}^+$ (the opening cannot equal $\pi$ because of (2.8)). The asymptotic expansion (3.14) can be determined then to higher orders in a straightforward way. (The equation for $G_k$ is (3.17), with r.h.s. a bilinear expression in $G_m$, $m < k$ (cf. also footnote 2)). The result is*

$$(3.19) \qquad g_C \sim G_0(s) + \frac{1}{x} G_1(s) + \ldots + \frac{1}{x^n} G_n(s) + \ldots$$

*with*

$$(3.20) \qquad G_0 = \frac{12s}{(s-1)^2}, \quad G_1 = -\frac{15s^3 + 175s^2 + 30s - 2}{10s(s-1)^3}, \quad \ldots$$

*It is also worth noting that for all $G_k(t)$ to be rational functions, the value $-392/625$ in (1.3) is crucial. A different numerical value produces a valid expansion (3.19) but with $G_k$, $k \geq 5$ not rational.*

**Remark 13.** *The arguments that were used to determine the asymptotic position of poles work even if $C$ is not large, as long as $|x|$ is large, yielding the location of poles in the neighborhood of the imaginary axis, for large $x$ (while on the positive real line $C$ small implies $x_p(C)$ small, by Proposition 1).*

## Appendix A. All solutions of P1 are meromorphic.

The fact that the only singularities of solutions of P1 are double poles was indicated by Painlevé, but it was not until the advent of the isomonodromic approach (see [1]) that a rigorous proof was given. For completeness, we present an independent proof of this crucial property of P1 (and thus of (1.3), for $x \neq 0$), based on the idea that near poles, the equation P1 is meromorphically equivalent to an elliptic equation [10].

**Notations:** $z_0$ denotes a point where initial conditions for $y$ are given, $\tilde{z}$ will denote singular points of $y$. Let $\tilde{D}$ be a disk without a discrete set of points. $\Lambda(z_1, z_2)$ will denote smooth curves in $\tilde{D}$, with the endpoints $z_1$ and $z_2 \in \text{Clos}(\tilde{D})$ removed, such that the length of any subsegment $\Lambda(z_3, z_4)$ of $\Lambda(z_1, z_2)$ is less than $2|z_3 - z_4|$.

**Lemma 14.** *In a neighborhood of any point $\tilde{z} \in \mathbb{C}$ there exists a one-parameter family of solutions $h$ of (1.1), having a double pole at $\tilde{z}$, of the form $h(z) = (z - \tilde{z})^{-2} A(z)$ where $A(z)$ is analytic at $\tilde{z}$, and $A(\tilde{z}) = 1$.*

*Proof.* Substituting $y = \sum_{k=-2}^{\infty} c_k (z - \tilde{z})^k$ in (1.1) we get

$$(A.1) \qquad (k+1)(k+2) c_{k+2} = 6 \sum_{m=-2}^{k+2} c_m c_{k-m} + \tilde{z} \delta_{k=0} + \delta_{k=1}$$

so that $c_{-2} = 0$ or $1$; choosing $c_{-2} = 1$ we have $c_{-1} = 0$ and

$$(A.2) \quad [(k+1)(k+2) - 12]c_{k+2} = 6\sum_{m=0}^{k} c_m c_{k-m} + \tilde{z}\delta_{k=0} + \delta_{k=1}$$

Carrying on the calculation we get

$$(A.3) \quad \left(c_{-2}, c_{-1}, c_0, c_1, c_2, c_3, c_4\right) = \left(1, 0, 0, 0, \frac{-1}{10}\tilde{z}, -\frac{1}{6}, c_4\right)$$

i.e., $c_4$ is a free constant (*Note*: The specific form of the equation (1.1) is crucial in obtaining an integer power series solution. Generic perturbations of (1.1) lead to an inconsistent equation for $c_4$.

Finally, for some $\rho > 0$ and all $k \geq 0$ we have $|c_k| \leq \frac{1}{3}(k+1)\rho^{k+2}$. Indeed, let $\rho$ be such that $|c_i| \leq \frac{1}{3}(i+1)\rho^{i+2}$ for $i \leq 2$. Then, assuming this property for $c_0, c_1, ..., c_{k+1}$ we have

$$|c_{k+2}| \leq \frac{1}{3}\frac{(k+1)(k+2)(k+3)}{(k-2)(k+5)}\rho^{k+4} \leq (k+3)\rho^{k+4}$$

□

**Lemma 15.** *If $y$ is a solution of (1.1), analytic at $z_0$, then the radius of analyticity at $z_0$ is at least $R = 1/\max\{|y(z_0)|^{1/2}, |y'(z_0)/2|^{1/3}, |y_0^2 + z_0/6|^{1/4}\}$.*

*Proof.* With $y = \sum_{k=0}^{\infty} c_k(z - z_0)^k$ we have

$$(A.4) \quad (k+1)(k+2)c_{k+2} = 6\sum_{m=0}^{k} c_m c_{k-m} + z_0\delta_{k=0} + \delta_{k=1}$$

By the definition of $R$ we have $|c_j| \leq (j+2)R^{-j-2}$ for $j \leq 2$. Inductively, if the inequality holds $\forall j \leq k+1$ we have

$$|c_{k+2}| \leq \frac{(k+1)(k+2)(k+3)R^{-k-4}}{(k+1)(k+2)} = (k+3)R^{-k-4}$$

□

**Lemma 16.** *Let $y$ be a solution of (1.1), analytic at any $z \in \Lambda(z_0, \tilde{z})$ and that $y(z) = \sum_{k=-2}^{4} a_k(z - \tilde{z})^k + O(z - \tilde{z})^5$ as $z \to \tilde{z}$ along $\Lambda(z_0, \tilde{z})$, with $a_{-2} \neq 0$.*
*Then $y$ is meromorphic at $\tilde{z}$, and $y$ is uniquely determined by $\tilde{z}$ and $a_4$.*

*Proof.* Since $\sum_{k=-2}^{4} a_k(z - \tilde{z})^k$ is an $O(z - \tilde{z})^5$ formal solution of P1, for $i \leq 3$ we must have $a_k = c_k$ defined in (A.3) (note: since $y$ solves (1.1), the asymptotic expansion is differentiable). Let $y_0(z)$ be the solution given in Lemma 14 with $c_4 = a_4$. Let $v(z) = y_0(z) - (z - \tilde{z})^{-2}$ and $f(z) = y(z) - y_0(z)$. Then $v = O((z - \tilde{z})^2), f = O(z - \tilde{z})^5$ and

$$(A.5) \quad f'' - \frac{12}{(z-\tilde{z})^2}f = 12fv + 6f^2$$

or

$$f(z) = C_1(z-\tilde{z})^{-3} + C_2(z-\tilde{z})^4 - \frac{7}{(z-\tilde{z})^3}\int_{\tilde{z}}^{z}(s-\tilde{z})^4 g(s)ds$$

(A.6)
$$+\frac{(z-\tilde{z})^4}{7}\int_{\tilde{z}}^{z}(s-\tilde{z})^{-3}g(s)ds = \mathcal{N}(f)$$

where $g = 12fv + 6f^2 = O(z-\tilde{z})^7$ near $\tilde{z}$. Thus $C_1 = C_2 = 0$. Let $z_1 \in \Lambda$ and $\mathcal{S}$ be the Banach space of continuous functions on $\Lambda(z_1, \tilde{z})$ with the norm $\sup_\Lambda |(z-\tilde{z})^{-5} f(z)|$. It is immediate to check that if $z_1 - \tilde{z}$ is small enough $\mathcal{N}$ is contractive in the unit ball $B \subset \mathcal{S}$. The solution of (A.6) is thus unique in $B$, $f = 0$. □

Let $y$ be a solution of (1.1), and $D$ an open disk where $y$ is meromorphic. By definition, the singularities of $y$ in $D$ are isolated. Let $\tilde{z} \in \partial D$ be a singular point of $y$. We will show that $\tilde{z}$ is a pole: *this implies that $y$ is meromorphic on $\mathbb{C}$.*

Let $\tilde{D}$ be $D$ without the poles of $y$.

**Remark 17.** *(i) As $z \to \tilde{z}$ along $\Lambda(z, \tilde{z})$ we have $|y(z)| + |y'(z)| \to \infty$. (ii) $y$ is unbounded as $z \to \tilde{z}$.*

*Proof.* (i) Since $y$ is singular at $\tilde{z}$ then by Lemma 15 $y(z)$ or $y'(z)$ must be large for $z$ near $\tilde{z}$. (ii) Multiplication of (1.1) by $y'$ and integration yields:

(A.7)
$$(y')^2 = 4y^3 + 2zy - 2\int_a^z y(s)ds + C$$

If $y$ were bounded along $\Lambda$, then by (A.7) $y'$ would be bounded as well. □.

**Proposition 18.** *Let $z_0 \in D$ be a regular point of $y$ and $\Lambda = \Lambda(z_0, \tilde{z}) \subset D$. We have $y^{-3}(s)\int_a^s y(t)dt \to 0$, $|y(s)| \to \infty$, and $|y'(s)| \to \infty$ as $s \to \tilde{z}$ on $\Lambda$.*

*Proof.* For *any* continuous function $f$ unbounded at $\tilde{z}$ we have

$$\liminf_{z\to\tilde{z}} |f|^{-3}\int_a^z |f|(t)d|t| \leq \liminf_{z\to\tilde{z}} |f|^{-3} \max_{\Lambda(a,z)} |f| \operatorname{length}(\Lambda) = 0$$

Thus, if we had $y^{-3}(s)\int_a^s y(t)dt \not\to 0$, then for any $\gamma > 0$ small enough

(A.8)
$$|\int_a^{s_n} y(t)dt| = \gamma|y^3(s_n)|$$

for some sequence $s_n \to \tilde{z}$ on $\Lambda$. If $\gamma < \frac{1}{2}$, by (A.8), (A.7) and Remark 17 we have

(A.9) $\quad |y(s_n)| \to \infty$ and $|y'(s_n)| \to \infty$ and $\int_{z_0}^{s_n} |y(t)|d|t| \to \infty$

Consider now the solution $h$ of the equation

(A.10)
$$(h')^2 = 4h^3 + 2s_n h + \tilde{C}$$

satisfying $h(0) = y(s_n), h'(0) = y'(s_n)$ (with $\tilde{C}$ chosen accordingly in (A.10)). It follows from (A.9) that $\tilde{C} = C - 2\int_a^{s_n} y(s)ds \to \infty$. Define the function $H$ by

$$\alpha H(\beta z) := h(z), \text{ where } \alpha = \tilde{C}^{1/3}, \ \beta = \tilde{C}^{1/6}$$

It follows that

$$(H')^2 = 4H^3 - 2\theta H + 1, \text{ where } \theta = \tilde{C}^{-2/3} \tag{A.11}$$

$H$ is thus an elliptic function with with periods $\Pi_{1,2} = \Pi + o(\theta)$ for some $\Pi \in \mathbb{C}$ (independent of $\theta$). By (A.8) $|y(s_n)| = (2\gamma^{-1/3})|\tilde{C}^{1/3}|(1+o(1))$, then using (A.7), we get $|y'(s_n)| = (4\gamma^{-1}e^{i\phi_n} + 2)^{-1/2}|\tilde{C}^{1/2}|(1+o(1))$ for some $\phi_n \in \mathbb{R}$.

Therefore $|H(0)| = (2\gamma)^{-1/3} + o(1)$, and $|H'(0)| = (4\gamma^{-1}e^{i\phi_n} + 2)^{-1/2}| + o(1)$. With the substitution $\alpha Y(\beta z) := y(s_n + z)$ we have

$$(Y')^2 = 4Y^3 + 2\theta Y + 1 + 2\epsilon\left(zY - \int_0^z Y(t)dt\right), \text{ where } \epsilon = \tilde{C}^{-5/6},$$
$$\text{or } Y'' = 6Y^2 + \theta + \epsilon s \tag{A.12}$$

We substitute $Y = H \circ G$; $G(z) = z + g(z)$ in (A.12) with $g(0) = g'(0) = 0$ and

$$2g'H'' \circ G + (g')^2 H'' \circ G + g''H' \circ G = \epsilon z$$
$$\text{or } 2g'(H' \circ G)' + g''H' \circ G = \epsilon z + (g')^2 H'' \circ G \tag{A.13}$$

which can be rewritten in the integral form ($g(z) = \int_{\Lambda(0,z)} w(t)dt, G = z + g$):

$$w = (H' \circ G)^{-2} \int_{\Lambda(0,z)} H'(G(u))\left[\epsilon u + H''(G(u))w^2(u)^2\right] du = \mathcal{N}(w) \tag{A.14}$$

Let $\Delta_\epsilon$ be the disk centered at 0, of diameter $|\epsilon|^{-1/4}$, and $\tilde{\Delta}_\epsilon \ni 0$ be $\Delta_\epsilon$ without the disks of diameter $d(\gamma)$ around the poles and zeros of $H'(z)$. Let $\Lambda \subset \Delta_\epsilon$, let $a \in (1/4, 1/2)$. Let $\mathcal{S}$ be the Banach space of continuous functions on $\Lambda$ in the sup norm and $B = \{w : \|w\| \le |\epsilon|^a\} \subset \mathcal{S}$. If $w \in B$, then $|g| \le |\int_{\Lambda(0,z)} w(t)dt| \le |\epsilon|^{a-1/4}$. $B$ is preserved by $\mathcal{N}$ provided that

$$|\epsilon|\, Const(\gamma)|\epsilon|^{-1/2} + Const(\gamma)|\epsilon|^{2a-1/4} \le |\epsilon|^a \tag{A.15}$$

which holds for small $\epsilon$. Similarly for $u, v \in B$, $\|\mathcal{N}(u) - \mathcal{N}(v)\| = O(\epsilon^{a-1/4})\|u - v\|$. $\mathcal{N}$ has thus a unique fixed point in $B$.

Therefore $|g(z)| \le |\epsilon|^{a-1/4}$ for $z \in \tilde{\Delta}_\epsilon$ and thus (when $\epsilon$ is small)

$$\sup_{\tilde{\Delta}_\epsilon} |\alpha^{-1}y(s_n + \beta^{-1}\cdot) - H(\cdot)| = \sup_{\tilde{\Delta}_\epsilon} |H(G(\cdot)) - H(\cdot)| = O(\epsilon^{a-1/4}) \tag{A.16}$$

Note now that as $s_n \to \tilde{z}$ we have $\beta^{-1} = \tilde{C}^{-1/6} \to 0$ while $\alpha \to \infty$ and $|\beta^{-1}\epsilon^{-1/4}| \propto \tilde{C}^{1/24} \to \infty$. Thus $L = s_n + \alpha^{-1}\tilde{\Delta}_\epsilon$ is an $\alpha^{-1}$-lattice for $D$ (i.e., $\forall z \in D, \exists z_1 \in L : |z - z_1| \leq \alpha^{-1}$), and (A.16) implies $y$ is unbounded near $z_0$, contradiction.

Since $y^{-3}(s) \int_a^s y(t)dt \to 0$, the fact that $\min\{|y|, |y'|\} \to \infty$ now follows from Remark 17 and (A.7). $\square$

**Proposition 19.** *$\tilde{z}$ is a second order pole.*

*Proof.* The substitution $y = 1/v^2$ transforms (A.7) into

$$(A.17) \qquad (v')^2 = 1 + \frac{1}{2}v(z)^4 z - \frac{1}{2}v(z)^6 \int_{z_0}^z v(s)^{-2} ds + Cv(z)^6$$

It follows from Proposition 18 that along a curve $\Lambda(z_0, \tilde{z})$ avoiding the poles of $y$ we have $v^6 \int_{z_0}^z v^{-2}(s)ds \to 0$. With $\zeta = z - \tilde{z}$, (A.17) implies $(v')^2 = 1 + o(1)$ so that $v^2 = \zeta^2(1+o(1))^2$ near $\tilde{z}$. Using (A.17) again we now get $v^2 = (\zeta + \sum_{k=2}^6 v_k \zeta^k + o(\zeta^6))$. Lemma 15 implies that $y$ is meromorphic at $\tilde{z}$. $\square$.

**Acknowledgments.** The author would like to thank Professors Percy Deift and Martin Kruskal for very interesting discussions on this subject.

OVIDIU COSTIN
DEPARTMENT OF MATHEMATICS, THE UNIVERSITY OF CHICAGO, 5734 S. UNIVERSITY AVE., CHICAGO IL 60637
TEMPORARY ADDRESS: MSRI, 1000 CENTENNIAL DR., BERKELEY, CA 94720
*E-mail address*: costin@math.uchicago.edu

---

[3]Preprints available at http://www.math.uchicago.edu/~costin